\documentclass[12pt,a4paper]{article}

\usepackage{amsmath}
\usepackage{amssymb}
\usepackage{amsthm}
\usepackage{t1enc}
\usepackage[latin1]{inputenc}
\usepackage[english]{babel}
\usepackage{amsfonts}
\usepackage{latexsym}
\usepackage{bm}
\usepackage{setspace}
\usepackage[toc,page]{appendix}
\usepackage{graphicx}
\usepackage{enumerate}
\usepackage[numbers,sort&compress]{natbib}
\usepackage{tikz}
\usepackage{csquotes}
\usepackage{doi}

\usepackage{subfigure}
\usepackage{floatrow}
\usepackage{mathtools}
\usepackage{caption}

\DeclarePairedDelimiter{\floor}{\lfloor}{\rfloor}

\newtheorem{theorem}{Theorem}
\newtheorem{lemma}{Lemma}
\newtheorem{corollary}{Corollary}
\newtheorem{proposition}{Proposition}

\newtheorem{construction}{Construction}

\usepackage{color}   
\usepackage{hyperref}
\hypersetup{
    colorlinks=true, 
    linktoc=all,     
    linkcolor=blue,  
    hidelinks,
}

\renewenvironment{proof}{{\noindent\bfseries Proof. }}

\usepackage{tikz}
\usetikzlibrary{shapes.geometric}
\tikzset{
    v/.style={draw,fill,circle,inner sep=1.5},
    e/.style={draw, thick}
}

\allowdisplaybreaks

\title{Solution to a problem on isolation of $3$-vertex paths}
\author{Karl Bartolo\footnote{Email address: karl.bartolo.16@um.edu.mt} \quad \quad Peter Borg\footnote{Email address: peter.borg@um.edu.mt} \quad \quad Dayle Scicluna\footnote{Email address: dayle.scicluna.09@um.edu.mt} \\[5mm]
{\normalsize Department of Mathematics} \\
{\normalsize Faculty of Science} \\
{\normalsize University of Malta}\\
{\normalsize Malta}
}
\date{}

\begin{document}

\maketitle

\begin{abstract}
The $3$-path isolation number of a connected $n$-vertex graph $G$, denoted by $\iota(G,P_3)$, is the size of a smallest subset $D$ of the vertex set of $G$ such that the closed neighbourhood $N[D]$ of $D$ intersects the vertex sets of the $3$-vertex paths of $G$, meaning that no two edges of $G-N[D]$ intersect. Zhang and Wu proved that $\iota(G,P_3) \leq 2n/7$ unless $G$ is a $3$-path or a $3$-cycle or a $6$-cycle. The bound is attained by infinitely many graphs having induced $6$-cycles. Huang, Zhang and Jin proved that if $G$ has no $6$-cycles, or $G$ has no induced $5$-cycles and no induced $6$-cycles, then $\iota(G, P_3) \leq n/4$ unless $G$ is a $3$-path or a $3$-cycle or a $7$-cycle or an $11$-cycle. They asked if the bound still holds asymptotically for connected graphs having no induced $6$-cycles. More precisely, taking $f(n)$ to be the maximum value of $\iota(G,P_3)$ over all connected $n$-vertex graphs $G$ having no induced $6$-cycles, their question is whether $\limsup_{n \to\infty}\frac{f(n)}{n} = \frac{1}{4}$. We verify this by proving that $f(n) = \left \lfloor (n+1)/4 \right \rfloor$. Therefore, induced $6$-cycles solely account for the jump to the bound $2n/7$. The proof hinges on further proving that if $G$ is such a graph and $\iota(G, P_3) = (n+1)/4$, then $\iota(G-v, P_3) < \iota(G, P_3)$ for each vertex $v$ of $G$. This new idea promises to be of further use. We also prove that if the maximum degree of such a graph $G$ is at least $5$, then $\iota(G,P_3) \leq n/4$. 
\end{abstract}

\section{Introduction}

Unless stated otherwise, we use small letters to denote non-negative integers or elements of sets, and capital letters to denote sets or graphs. The set of positive integers is denoted by $\mathbb{N}$. For $n \in \{0\} \cup \mathbb{N}$, $[n]$ denotes the set $\{i\in \mathbb{N} \colon i \leq n\}$. Note that $[0]$ is the empty set $\emptyset$. Arbitrary sets are taken to be finite. For a set $X$, $\binom{X}{2}$ denotes the set of $2$-element subsets of $X$ (that is, $\binom{X}{2} = \{ \{x,y \} \colon x,y \in X, x \neq y \}$). A $2$-element set $\{x,y\}$ may be represented by $xy$. For standard terminology in graph theory, the reader is referred to~\cite{West}. The terminology, notation and motivation for this paper are largely provided in our related paper~\cite{BBS}, but we provide them here in full for the convenience of the reader who has not read~\cite{BBS}.

Every graph $G$ is taken to be \emph{simple}, that is, $G$ is a pair $(V(G), E(G))$ such that the \emph{vertex set} $V(G)$ and \emph{edge set} $E(G)$ satisfy $E(G) \subseteq \binom{V(G)}{2}$. If $|V(G)| = n$, then $G$ is called an \emph{$n$-vertex graph}. For $v \in V(G)$, $N_{G}(v)$ denotes the set of neighbours of $v$ in $G$, $N_{G}[v]$ denotes the closed neighbourhood $N_{G}(v) \cup \{ v \}$ of $v$, and $d_{G}(v)$ denotes the degree $|N_{G} (v)|$ of $v$. The maximum degree $\max\{d_G(v) \colon v \in V(G)\}$ of $G$ is denoted by $\Delta(G)$. If $\Delta(G) \leq 3$, then $G$ is said to be \emph{subcubic}. For a subset $X$ of $V(G)$, $N_G[X]$ denotes the closed neighbourhood $\bigcup_{v \in X} N_G[v]$ of $X$, $N_G(X)$ denotes the open neighbourhood $N_G[X]\setminus X$ of $X$, $G[X]$ denotes the subgraph of $G$ induced by $X$, and $G - X$ denotes the graph obtained by deleting $X$ from $G$. Thus, $G[X] = (X,E(G) \cap \binom{X}{2})$ and $G - X = G[V(G) \setminus X]$. We may abbreviate $G-\{x\}$ to $G-x$. Where confusion does not arise, the subscript $G$ may be omitted.

Consider two graphs $G$ and $H$. If $G$ is a copy of $H$, then we say that $G$ is an \emph{$H$-copy}, and write $G \simeq H$. If $V(H) \subseteq V(G)$ and $E(H) \subseteq E(G)$, then we say that $G$ \emph{contains} $H$ and that $H$ is a \emph{subgraph of $G$}. If $H = G[X]$ for some $X \subseteq V(G)$, then $H$ is called an \emph{induced subgraph of $G$}.
  
For $n \geq 1$, $K_n$ and $P_n$ denote the graphs $([n], \binom{[n]}{2})$ and $([n], \{\{i,i+1\} \colon i \in [n-1]\})$, respectively. For $n \geq 3$, $C_n$ denotes $([n],E(P_n) \cup \{\{n,1\}\})$. A $K_n$-copy is called an \emph{$n$-clique}, a $P_n$-copy is called an \emph{$n$-path} or simply a \emph{path}, and a $C_n$-copy is called an \emph{$n$-cycle} or simply a \emph{cycle}. 

If $\mathcal{F}$ is a set of graphs and $F$ is a copy of a graph in $\mathcal{F}$, then we call $F$ an \emph{$\mathcal{F}$-graph}. A subset $D$ of $V(G)$ is called an \emph{$\mathcal{F}$-isolating set of $G$} if $N[D]$ intersects the vertex sets of the $\mathcal{F}$-graphs contained by $G$. Thus, $D$ is an $\mathcal{F}$-isolating set of $G$ if and only if $G-N[D]$ contains no $\mathcal{F}$-graph. The size of a smallest $\mathcal{F}$-isolating set of $G$ is called the \emph{$\mathcal{F}$-isolation number of $G$} and denoted by $\iota(G, \mathcal{F})$. If $\mathcal{F} = \{F\}$, then we may replace $\mathcal{F}$ in the aforementioned terms and notation by $F$.

The study of $\mathcal{F}$-isolating sets, initiated by Caro and Hansberg \cite{CaHa17}, vastly widens the popular study of dominating sets~\cite{C, CH, HHS, HHS2, HL, HL2} in a natural way. Indeed, $D$ is a \emph{dominating set of $G$} (that is, $D \subseteq V(G) = N[D]$) if and only if $D$ is a $K_1$-isolating set of $G$, so the \emph{domination number of $G$} (the size of a smallest dominating set of $G$), denoted by $\gamma(G)$, is $\iota(G,K_1)$.

Let $G$ be a connected $n$-vertex graph. A classical result of Ore~\cite{Ore} is that $\gamma(G) \leq n/2$ if $G \not\simeq K_1$ (see~\cite{HHS}). This is one of the earliest domination results. While deleting the closed neighbourhood of a dominating set produces the graph with no vertices, deleting the closed neighbourhood of a $K_2$-isolating set produces a graph with no edges. In \cite{CaHa17}, Caro and Hansberg proved that $\iota(G, K_2) \leq n/3$ unless $G \simeq K_2$ or $G \simeq C_5$ (this was independently proved by \.{Z}yli\'{n}ski~\cite{Z}; see \cite{LMS, BG} for the cases of equality). The same paper featured several problems, some of which have been settled by the present second author and others \cite{Borg1,Borgrsc,BFK} (see also \cite{BFK2, Borgrc2}). Domination and isolation have been particularly investigated for maximal outerplanar graphs~\cite{BK, BK2, CaWa13, CaHa17, Ch75, DoHaJo16, DoHaJo17, HeKa18, LeZuZy17, Li16, MaTa96, To13, KaJi}, mostly due to connections with Chv\'{a}tal's Art Gallery Theorem~\cite{Ch75}. As in the development of domination, isolation is expanding in various directions; for example, total isolation and an isolation game have been treated in \cite{BGH} and \cite{BDJKR}, respectively.

Again consider a connected $n$-vertex graph $G$. With the establishment of a sharp upper bound on each of $\iota(G,K_1)$ and $\iota(G,K_2)$, the natural subsequent target was a sharp upper bound on $\iota(G,F)$ for each of the two cases where $F$ is a connected $3$-vertex graph, or rather, $F \in \{K_3, P_3\}$. 

The sharp inequality $\iota(G,K_3) \leq n/4$ with $G \not\simeq K_3$ was established in \cite{Borg1} as a bound on the cycle isolation number, and is also given in \cite{BFK} as the special case $k = 3$ of the sharp inequality $\iota(G,K_k) \leq n/(k+1)$ with $G \not\simeq K_k$, in \cite{Borgcon} as an immediate consequence of the inequality $\iota(G,\mathcal{E}_3) \leq n/4$ with $\mathcal{E}_3 = \{H \colon H \mbox{ is a connected graph, } |E(H)| \geq 3\}$ and $G$ not being a $\{K_3, C_7\}$-graph, and in \cite{Borgrsc} as an immediate consequence of a result on isolation of regular graphs, stars and $k$-chromatic graphs. 

The significance of a $P_3$-isolating set $D$ is that $E(G-N[D])$ is a \emph{matching}, meaning that no two edges of $G-N[D]$ have a common vertex (and hence $\Delta(G-N[D]) \leq 1$). If $G$ is not a $\{K_3, P_3, C_6\}$-graph, then $\iota(G,P_3) \leq 2n/7$. This was established by Zhang and Wu in \cite{ZW}, and independently by the present second author in \cite{Borgcon}, which actually established the refined inequality $\iota(G,P_3) \leq (4n-r)/14$ with $r$ being the number of leaves (vertices of degree $1$) of $G$, and $G$ not being a copy of one of six particular graphs. The graphs that attain the bound $2n/7$ were determined independently in \cite{CLWX} and \cite{CZZ}. Infinitely many of them have induced $6$-cycles. Improving another result of Zhang and Wu in \cite{ZW}, Huang, Zhang and Jin \cite{HZJ} showed that $\iota(G,P_3) \leq n/4$ if $G$ is not a $\{P_3,C_3,C_7,C_{11}\}$-graph and $G$ either contains no $6$-cycles or contains $6$-cycles but not induced $5$-cycles or induced $6$-cycles. They asked if the bound still holds asymptotically for connected graphs having no induced $6$-cycles. More precisely, let 
\[f(n) = \max\{\iota(G,P_3) \colon G \mbox{ is a connected $n$-vertex graph that has no induced $6$-cycles}\}.\]
The question in \cite[Problem~5.1]{HZJ} is whether 
$$\limsup_{n \to\infty} \frac{f(n)}{n} = \frac{1}{4}.$$ 
Thus, the problem essentially is whether induced $6$-cycles solely account for the jump from the bound $n/4$ to the bound $2n/7$ on $\iota(G,P_3)$. In \cite{BBS}, we solved this problem for subcubic graphs, which need to be treated differently from other graphs. Throughout the rest of this section, it is to be assumed that $G$ has no induced $6$-cycles. We showed that if $G$ is subcubic, then $\iota(G,P_3) \leq n/4$ unless $G$ is a copy of one of $12$ particular graphs whose $P_3$-isolation number is $(n+1)/4$. In this paper, we solve the problem completely. Further investigation led us to the discovery of many infinite sets of non-subcubic graphs $G$ with $\iota(G,P_3) = (n+1)/4$, one of which is given in the following construction (see Figure~\ref{FigExtremal}).

\begin{figure}[htb!]
    {\caption{The graph $B_4^*$ defined in Construction~\ref{ConsB}.}\label{FigExtremal}
    \centering
    \begin{tikzpicture}
        \node[v,label={east:$w_4$}] (1) at (4,-4.5) {};
        \node[v,label={north:$v_4$}] (2) at (3,-4) {};
        \node[v,label={north:$w_3$}] (3) at (2,-4) {};
        \node[v,label={north:$v_3$}] (4) at (1,-4) {};
        \node[v,label={south:$v_3'$}] (5) at (1,-5) {}; 
        \node[v,label={south:$w_3'$}] (6) at (2,-5) {};     
        \node[v,label={south:$v_4'$}] (7) at (3,-5) {} ;
        \node[v,label={north:$w_2$}] (8) at (0,-4) {};
        \node[v,label={north:$v_2$}] (9) at (-1,-4) {};
        \node[v,label={south:$w_2'$}] (10) at (0,-5) {};
        \node[v,label={south:$v_2'$}] (11) at (-1,-5){};         
        \node[v,label={north:$w_1$}] (12) at (-2,-4) {};
        \node[v,label={north:$v_1$}] (13) at (-3,-4) {} ;
        \node[v,label={south:$w_1'$}] (14) at (-2,-5) {};         
        \node[v,label={south:$v_1'$}] (15) at (-3,-5) {};  
            
        \draw[e] (1) edge (2);
        \draw[e] (1) edge (7);
        \draw[e] (2) edge (3);
        \draw[e] (2) edge (7);
        \draw[e] (3) edge (5);
        \draw[e] (3) edge (4); 
        \draw[e] (4) edge (5);
        \draw[e] (4) edge (6);
        \draw[e] (5) edge (6);  
        \draw[e] (6) edge (7);
        \draw[e] (4) edge (8);
        \draw[e] (5) edge (10);
        \draw[e] (10) edge (11);
        \draw[e] (8) edge (9); 
        \draw[e] (9) edge (11);
        \draw[e] (8) edge (11);
        \draw[e] (9) edge (10);
        		
        \draw[e] (9) edge (12);
        \draw[e] (11) edge (14);
        \draw[e] (12) edge (13);
        \draw[e] (13) edge (15); 
        \draw[e] (14) edge (15);
        \draw[e] (12) edge (15);
        \draw[e] (13) edge (14);
        		
    \end{tikzpicture}    
    }
\end{figure}

Let $K_4^-$ denote the graph with vertex set $V(K_4)$ and edge set $E(K_4) \setminus \{\{3,4\}\}$.

\begin{construction}\label{ConsB} \emph{Consider any $k \in \mathbb{N}$. Let $F_1, \dots, F_k$ be pairwise vertex-disjoint graphs, where $F_k \simeq K_3$ and $F_i \simeq K_4^-$ for each $i \in [k-1]$. For each $i \in [k-1]$, let $v_i, v_i', w_i, w_i'$ be the vertices of $F_i$ corresponding to the vertices $1, 2, 3, 4$ of $K_4^-$, respectively. Let $v_k, v_k', w_k$ be the vertices of $F_k$ corresponding to the vertices $1, 2, 3$ of $K_3$, respectively. Let $B_k^*$ be the connected $(4k-1)$-vertex graph with vertex set $\bigcup_{i=1}^{k} V(F_i)$ and edge set $\left( \bigcup_{i=1}^{k} E(F_i) \right) \cup \bigcup_{i=1}^{k-1} \{w_iv_{i+1}, w_i'v_{i+1}'\}$.}
\end{construction}

Note that $B_1^* \simeq K_3$ and $B_2^*$ is the graph $G_{7,6}$ in~\cite{BBS}. The following is one of the results proved in the next section. 

\begin{proposition}\label{PropB} For any $k \geq 1$, $B_k^*$ has no induced $6$-cycles and
\[\iota(B_k^*,P_3) = k = \frac{|V(B_k^*)|+1}{4}.\]
\end{proposition}
\noindent
Using Construction~\ref{ConsB}, for any $n \geq 1$, we now construct a connected $n$-vertex graph $B_n$ with no induced $6$-cycles and with $\iota(B_n,P_3) = \lfloor (n+1)/4 \rfloor$ (see Corollary~\ref{PropBcor}).

\begin{construction}\label{ConsB2} \emph{Consider any $n \in \mathbb{N}$. Let $k = \lfloor (n+1)/4 \rfloor$. Thus, $n+1 = 4k + r$ for some $r \in \{0, 1, 2, 3\}$. If $k = 0$, then let $B_n = P_n$. If $r = 0$, then let $B_n = B_k^*$. If $k \geq 1$ and $r \geq 1$, then let $F^* \simeq P_{r+1}$ such that $V(B_k^*) \cap V(F^*) = \{v_1\}$, and let $B_n$ be the connected $n$-vertex graph with vertex set $V(B_k^*) \cup V(F^*)$ and edge set $E(B_k^*) \cup E(F^*)$.}
\end{construction}
\noindent
Note that $\Delta(B_n) = 4$ for $n \geq 8$. Thus, there are infinitely many graphs $G$ with $\Delta(G) = 4$ and $\iota(G,P_3) > n/4$, and as mentioned above, there are 12 graphs $G$ with $\Delta(G) \leq 3$ and $\iota(G,P_3) > n/4$. We show that $\iota(G,P_3) \leq n/4$ if $\Delta(G) \geq 5$. In light of this fact, for any $h \geq 5$ and $n \geq 4(h-2)$, we construct a connected $n$-vertex graph $B_{n,K_3,h}$ with no induced $6$-cycles and with $\iota(B_{n,K_3,h}, P_3) = \lfloor n/4 \rfloor$ and $\Delta(B_{n,K_3,h}) = h$ (see Figure~\ref{FigB3} and Proposition~\ref{PropB3}). We modify the construction of the graph $B_{n,F}$ with $F = K_3$ in \cite{Borg1}.

\begin{figure}[htb!]
    {\caption{The graph $B_{30,K_3,6}$ defined in Construction~\ref{ConsB3}.}\label{FigB3}
    \centering
\begin{tikzpicture}
	\foreach \i in {0,...,6}{
    \node[v] (v\i1)  at (0.5+1.75*\i,1.5) {};
	\node[v] (v\i2)  at (1+1.75*\i,0) {};
    \node[v] (v\i3)  at (1.75*\i,0) {};
	\node[v] (v\i4)  at (0.5+1.75*\i,0.5) {};
		
	\draw[e] (v\i1) -- (v\i2);
	\draw[e] (v\i1) -- (v\i3);
	\draw[e] (v\i1) -- (v\i4);
	\draw[e] (v\i2) -- (v\i3);
	\draw[e] (v\i2) -- (v\i4);
	\draw[e] (v\i3) -- (v\i4);
	}
    
 	\draw[e] (v01) -- (v11);
 	\draw[e] (v01) .. controls (5/3,1.75) and (8.5/3,1.75)  .. (v21);
 	\draw[e] (v01) .. controls (2.1875,2.125) and (3.9375,2.125)  .. (v31);
	
	\draw[e] (v31) -- (v41);
	\draw[e] (v41) -- (v51);
	\draw[e] (v51) -- (v61);

	\node[v] (b1) at (10.125+1.75,1.5) {};
	\node[v] (b2) at (11+1.75,1.5) {};

	\draw[e] (v61) -- (b1);
	\draw[e] (b1) -- (b2);

\end{tikzpicture}
}
\end{figure}

\begin{construction}\label{ConsB3} \emph{Consider any $n, h \in \mathbb{N}$ with $h \geq 5$ and $n \geq 4(h-2)$. Let $a_n = \lfloor n/4 \rfloor$, let $b_{n} = n - 3a_{n}$ (so $a_{n} \leq b_{n} \leq a_{n} + 3$), let $F_1, \dots, F_{a_{n}}$ be copies of $K_3$ such that $P_{b_{n}}, F_1, \dots, F_{a_{n}}$ are pairwise vertex-disjoint, and let $B_{n,K_3,h}$ be the connected $n$-vertex graph with vertex set $[b_{n}] \cup \bigcup_{i=1}^{a_{n}} V(F_i)$ and edge set $(E(P_{b_{n}}) \setminus E(P_{h-2})) \cup \{\{1,i+1\} \colon i \in [h-3]\} \cup \bigcup_{i=1}^{a_{n}} (E(F_i) \cup \{\{i,v\} \colon v \in V(F_i)\})$.}
\end{construction}

We now state our main result, which is proved in Section~\ref{Sectionpf}.

\begin{theorem}\label{ThmP3beyond} Let $G$ be a connected $n$-vertex graph that has no induced $6$-cycles. Then, 
\begin{equation} \iota(G,P_3) \leq \left \lfloor \frac{n+1}{4} \right \rfloor, \label{thmfact1} 
\end{equation}
and equality holds if $G = B_n$. Moreover, the following statements hold:\medskip 
\\
(a) If $\iota(G,P_3) = \frac{n+1}{4}$, then for each $v \in V(G)$,
\[\iota(G-v,P_3) \leq \iota(G,P_3) - 1.\]
(b) If $h = \Delta(G) \geq 5$ or $\iota(G,P_3) \neq \frac{n+1}{4}$, then 
\begin{equation} \iota(G,P_3) \leq \left \lfloor \frac{n}{4} \right \rfloor, \label{thmfact2} 
\end{equation}
and equality holds if $n \geq 4(h-2)$ and $G = B_{n,K_3,h}$.
\end{theorem}
\noindent The proof uses ideas from \cite{Borgcon, Borgisdom}. The proof of (\ref{thmfact1}) and (\ref{thmfact2}) is based on part (a) of the theorem. This is a new idea that promises to be of further use. We mention that equality holds in (a) (by Lemma~\ref{LemIsol1}), but the inequality suffices for our purpose.

Theorem~\ref{ThmP3beyond} has the following immediate consequence.
\begin{corollary} For any $n \geq 1$,
\[f(n) = \left \lfloor \frac{n+1}{4} \right \rfloor.\]
\end{corollary}

It is worth pointing out that the bound in (\ref{thmfact1}) appears to be the first of its kind because, as indicated above, sharp upper bounds on $\iota(G,\mathcal{F})$ in terms of $n$ are typically of the form $\alpha n$ for some positive constant $\alpha < 1$. Despite this intriguing fact, we provide an affirmative answer to the above-mentioned question of Huang, Zhang and~Jin.

\begin{corollary}[Solution to {\cite[Problem 5.1]{HZJ}}]\label{CorLimsup} 
\[\limsup_{n \to\infty} \frac{f(n)}{n} = \frac{1}{4}.\] 
\end{corollary}
\proof{By Theorem~\ref{ThmP3beyond}, for any $n \in \mathbb{N}$, we have $\frac{f(n)}{n} \leq \frac{1}{4}$ if $4$ does not divide $n+1$, and $\frac{f(n)}{n} = \frac{1}{4} + \frac{1}{4n}$ if $4$ divides $n+1$. Let $n^+$ be the smallest integer in $\{k \in \mathbb{N} \colon k \geq n, \, 4 \mbox{ divides } k+1\}$. Thus, $\lim_{n \rightarrow \infty} \sup \left\{\frac{f(k)}{k} \colon k \geq n \right\}  = \lim_{n \rightarrow \infty} \left( \frac{1}{4} + \frac{1}{4n^+} \right)  = \frac{1}{4}$. \qed}

\section{Proof of Theorem~\ref{ThmP3beyond}}\label{Sectionpf}

We now prove Theorem~\ref{ThmP3beyond}. We often use the following two lemmas.

\begin{lemma}[\cite{Borg1}]\label{LemIsol1} If $G$ is a graph, $\mathcal{F}$ is a set of graphs, $X\subseteq V(G)$, and $Y\subseteq N[X]$, then \[\iota(G,\mathcal{F}) \leq |X| + \iota(G-Y,\mathcal{F}).\]
\end{lemma}

A \emph{component} of a graph $G$ is a maximal connected subgraph of $G$. Clearly, the components of $G$ are pairwise vertex-disjoint (that is, no two have a common vertex), and their union is $G$. The set of components of $G$ is denoted by $\mathrm{C}(G)$.  

\begin{lemma}[\cite{Borg1, Borgrsc}]\label{LemIsol2} If $G$ is a graph and $\mathcal{F}$ is a set of connected graphs, then \[\iota(G,\mathcal{F}) = \sum_{H \in \mathrm{C}(G)}\iota(H,\mathcal{F}).\]
\end{lemma}

\noindent \textbf{Proof of Proposition~\ref{PropB}.} We use induction on $k$. The result is trivial if $k = 1$. Suppose $k \geq 2$. Let $G = B_k^*$ and $n = |V(G)|$. We have $k = (n+1)/4$. 

Let $G^* = G - N[v_1]$. Then, $G^* = G - \{v_1, v_1', w_1, w_1'\} \simeq B_{k-1}^*$. By the induction hypothesis, $G^*$ has no induced $6$-cycles. Suppose $G[X] \simeq C_6$ for some $X \subseteq V(G)$. Then, $X \cap N[v_1] \neq \emptyset$ and $X \cap V(G^*) \neq \emptyset$. It clearly follows that $w_1v_2, w_1'v_2' \in E(G[X])$, so $X = \{w_1, v_2, w_1', v_2', x_1, x_2\}$ for some $x_1 \in \{v_1, v_1'\}$ and $x_2 \in \{v_1, v_1', w_2, w_2'\} \setminus \{x_1\}$. By symmetry, we may assume that $x_1 = v_1$. Since $G[\{v_1, v_1', w_1, w_1', v_2, v_2'\}] \not\simeq C_6$, we have $x_2 \neq v_1'$, so $x_2 \in \{w_2, w_2'\}$. We have $|E(G[X])| = 7$, contradicting $G[X] \simeq C_6$. Therefore, $G$ has no induced $6$-cycles.

Clearly, $\{v_i \colon i \in [k]\}$ is a $P_3$-isolating set of $B_k^*$, so $\iota(B_k^*, P_3) \leq k$. Let $D$ be a $P_3$-isolating set of $G$. The result follows if we show that $|D| \geq k$. 
Let $G' = G - N[v_k']$. We have $N(v_k') = \{w_{k-1}', v_k, w_k\}$ and $G' = B_{k-1}^*$.

\noindent\textbf{Case 1:} $v_k \in D$ or $v_k' \in D$. By symmetry, we may assume that $v_k' \in D$. Let $D^* =  D \cap N[v_k']$. Then, $v_k' \in D^*$. If $D^* = \{v_k'\}$, then $N[D^*] \cap V(G') = \emptyset$, and let $D' = D \setminus D^*$. If $D^* \neq \{v_k'\}$, then $|D^*| \geq 2$ and $N[D^*] \cap V(G') \subseteq N[v_{k-1}]$, and let $D' = (D \setminus D^*) \cup \{v_{k-1}\}$. Thus, $D'$ is a $P_3$-isolating set of $G'$, and $|D| \geq |D'| + 1$. By the induction hypothesis, $|D'| \geq k-1$, so $|D| \geq k$. 

\noindent\textbf{Case 2:} $v_k, v_{k-1}' \notin D$.

\noindent\textbf{Case 2.1:} $w_k \in D$. Let $D' = (D \setminus \{w_k\}) \cup \{v_k'\}$. Since $N[w_k] \subseteq N[v_k']$, $D'$ is a $P_3$-isolating set of $G$. We have $|D| = |D'| \geq k$ by Case~1.

\noindent\textbf{Case 2.2:} $w_k \notin D$. Since $D$ is a $P_3$-isolating set of $G$ and $v_k, v_k', w_k \notin D$, we have $w_{k-1} \in D$ or $w_{k-1}' \in D$. We may assume that $w_{k-1}' \in D$.  Since $G[\{w_{k-1}, v_k, w_k\}]$ is a copy of $P_3$ and has no vertex in $N[w_{k-1}']$, $D$ has some vertex $x$ in $\{v_{k-1}, v_{k-1}', w_{k-1}\}$. Let $D' = D \setminus \{w_{k-1}'\}$. Since $N[w_{k-1}'] \cap V(G') \subseteq N[x]$, $D'$ is a $P_3$-isolating set of $G'$. By the induction hypothesis, $|D'| \geq k-1$. We have $|D| = |D'| + 1 \geq k$. \qed{}

\begin{corollary}\label{PropBcor} For any $n \geq 1$, $B_n$ has no induced $6$-cycles and
\[\iota(B_n, P_3) = \left\lfloor \frac{n+1}{4} \right\rfloor.\]
\end{corollary}
\noindent
\textbf{Proof.} The result is trivial if $k = 0$. The result is Proposition~\ref{PropB} if $r = 0$. Suppose $k \geq 1$ and $r \geq 1$. It is immediate from Proposition~\ref{PropB} that $B_n$ has no induced $6$-cycles. Let $D$ be a smallest $P_3$-isolating set of $B_n$. Clearly, $\{v_i \colon i \in [k]\}$ is a $P_3$-isolating set of $B_n$, so $|D| \leq k$. Let $X = V(B_n) \setminus V(B_k^*)$. Thus, $X = V(F^*) \setminus \{v_1\}$. If $D \cap X = \emptyset$, then let $D' = D$. If $D \cap X \neq \emptyset$, then let $D' = (D \setminus X) \cup \{v_1\}$. Thus, $D'$ is a $P_3$-isolating set of $B_k^*$, and we have $k \geq |D| \geq |D'| \geq k$ by Proposition~\ref{PropB}. Therefore, $|D| = k$. \qed{}

\begin{proposition}\label{PropB3} For any $h \geq 5$ and $n \geq 4(h-2)$, $B_{n,K_3,h}$ has no induced $6$-cycles, $\Delta(B_{n,K_3,h}) = h$ and
\[\iota(B_{n,K_3,h}, P_3) = \left\lfloor \frac{n}{4} \right\rfloor.\]
\end{proposition}
\noindent
\textbf{Proof.} Let $G = B_{n,K_3,h}$. It is immediate from Construction~\ref{ConsB3} that $G$ has no induced $6$-cycles. We have $d(1) = |V(F_1)| + h-3 = h$ and $d(v) \leq |V(K_3)| + 2 = 5$ for each $v \in V(G) \setminus \{1\}$, so $\Delta(G) = h$ as $h \geq 5$. Note that $[a_n]$ is a $P_3$-isolating set of $G$, so $\iota(G, P_3) \leq a_n$. For each $i \in [a_n]$, $N(V(F_i)) = \{i\}$. Thus, if $D$ is a $P_3$-isolating set of $G$, then $D \cap (V(F_i) \cup \{i\}) \neq \emptyset$ for each $i \in [a_n]$. Therefore, $\iota(G, P_3) = a_n$. \qed{}
\\
\\
\textbf{Proof of Theorem~\ref{ThmP3beyond}.}
By Corollary~\ref{PropBcor}, equality in (\ref{thmfact1}) holds if $G = B_n$. By Proposition~\ref{PropB3}, equality in (\ref{thmfact2}) holds if $h = \Delta(G) \geq 5$, $n \geq 4(h-2)$ and $G = B_{n,K_3,h}$. 

We now prove (\ref{thmfact1}), (\ref{thmfact2}) and (a), using induction on $n$. Note that since $G$ has no induced $6$-cycles, for any $X \subseteq V(G)$, $G-X$ has no induced $6$-cycles. Also note that (\ref{thmfact1}) implies that if $\iota(G, P_3) \neq (n+1)/4$, then $4 \iota(G, P_3) < n+1$, so $\iota(G, P_3) \leq n/4$.

The result is trivial if $\Delta(G) \leq 1$. Suppose $\Delta(G) \geq 2$. Then, $n \geq 3$. If $n=3$, then $\iota(G,P_3) = 1 = (n+1)/4$ and $\iota(G-v,P_3) = 0$ for any $v \in V(G)$. Suppose $n \geq 4$. 

We first prove (a). Thus, suppose $\iota(G,P_3) = (n+1)/4$. 
Let $\mathcal{I} = \mathrm{C}(G-v)$, $\mathcal{I}' = \{I \in \mathcal{I} \colon \iota(I, P_3) = (|V(I)|+1)/4\}$ and $\mathcal{I}^*=\mathcal{I} \setminus \mathcal{I}'$. By the induction hypothesis, $\iota(I, P_3) \leq (|V(I)|+1)/4$ for each $I \in \mathcal{I}$. For each $I \in \mathcal{I}^*$, we have $\iota(I, P_3) < (|V(I)|+1)/4$, so $4\iota(I, P_3) < |V(I)|+1$, and hence $\iota(I, P_3) \leq |V(I)|/4$. Since $G$ is connected, for each $I \in \mathcal{I}$, $vx_I \in E(G)$ for some $x_I \in V(I)$. Suppose $\mathcal{I}' \neq \emptyset$. By the induction hypothesis, for each $I \in \mathcal{I}'$, $\iota(I - x_I, P_3) \leq \iota(I, P_3) - 1$, so $\iota(I - x_I, P_3) \leq (|V(I)|-3)/4$. Let $Y = \{v\} \cup \{x_I \colon I \in \mathcal{I}'\}$. Then, $Y \subseteq N[v]$ and $\mathrm{C}(G-Y) = \mathcal{I}^* \cup \bigcup_{I \in \mathcal{I}'} \mathrm{C}(I - x_I)$. By Lemma~\ref{LemIsol1} (with $X = \{v\}$) and Lemma~\ref{LemIsol2}, we have
\begin{align*} \iota(G, P_3) &\leq 1 + \iota(G-Y,P_3) = 1 + \sum_{I \in \mathcal{I}^*} \iota(I, P_3) + \sum_{I \in \mathcal{I}'} \sum_{H \in \mathrm{C}(I - x_I)} \iota(H, P_3) \\
&= 1 + \sum_{I \in \mathcal{I}^*} \iota(I, P_3) + \sum_{I \in \mathcal{I}'} \iota(I - x_I, P_3) \leq \frac{4}{4} + \sum_{I \in \mathcal{I}^*} \frac{|V(I)|}{4} + \sum_{I \in \mathcal{I}'} \frac{|V(I)|-3}{4} \\
&=\frac{4-3|\mathcal{I}'|+n-1}{4}=\frac{n-3(|\mathcal{I}'|-1)}{4}<\frac{n+1}{4},
\end{align*}
which contradicts $\iota(G,P_3) = (n+1)/4$. Thus, $\mathcal{I}' = \emptyset$, and hence $\mathcal{I} = \mathcal{I}^*$. By Lemma~\ref{LemIsol2},
\[\iota(G-v, P_3) \leq \sum_{I \in \mathcal{I}} \frac{|V(I)|}{4} = \frac{n-1}{4} < \frac{n+1}{4}= \iota(G,P_3).\]
Since $\iota(G-v, P_3)$ and $\iota(G, P_3)$ are integers, $\iota(G-v, P_3) \leq \iota(G, P_3) - 1$.

We now prove (\ref{thmfact1}) and (\ref{thmfact2}). Since $\iota(G,P_3)$ is an integer, it suffices to show that $\iota(G,P_3) \leq (n+1)/4$ if $\Delta(G) \leq 4$, and that $\iota(G,P_3) \leq n/4$ if $\Delta(G) \geq 5$. 

Suppose $\Delta(G) = 2$. Since $G$ is connected, we have $G \simeq P_n$ or $G \simeq C_n$. Recall that $n \geq 4$. If $G = P_n$, then $\{4k \colon k\in [\floor{n/4}]\}$ is a $P_3$-isolating set of $G$. If $G = C_n$, then $n\neq 6$ (as $G$ has no induced $6$-cycles) and $\{5k-4 \colon k \in \mathbb{N}, \, 5k-4 \leq n\}$ is a $P_3$-isolating set of size $\floor{(n+4)/5}$. Clearly, $\floor{(n+4)/5} \leq (n+1)/4$ if $n \in \{4,5,7,8,9,10\}$, and $(n+4)/5 \leq (n+1)/4$ if $n\geq 11$. Therefore, $\iota(G,P_3) \leq (n+1)/4$. 

Now suppose $\Delta(G) \geq 3$. Let $v\in V(G)$ with $d(v) = \Delta(G)$. If $N[v]=V(G)$, then $\iota(G,P_3)=1 \leq n/4$. Suppose $N[v] \neq V(G)$. Let $G'=G-N[v]$. Since $\Delta(G) \geq 3$ and $V(G') \neq \emptyset$, $n\geq 5$. 

Let $\mathcal{H} = \mathrm{C}(G')$, $\mathcal{H}' = \{H \in \mathcal{H} \colon \iota(H,P_3) = (|V(H)|+1)/4\}$ and $\mathcal{H}^*=\mathcal{H}\setminus\mathcal{H}'$. Thus, for each $H \in \mathcal{H}'$, $|V(H)| = 4k+3$ for some $k \in \{0\} \cup \mathbb{N}$. By the induction hypothesis, $\iota(H,P_3) \leq (|V(H)|+1)/4$ for each $H \in \mathcal{H}$. By the argument above for $\mathcal{I}^*$, $\iota(H,P_3)\leq |V(H)|/4$ for each $H \in \mathcal{H}^*$.

If $\mathcal{H}' = \emptyset$, then by Lemmas~\ref{LemIsol1} and~\ref{LemIsol2},
\begin{align*}
\iota(G,P_3) &\leq 1 + \iota(G',P_3)\leq \frac{|N[v]|}{4} + \sum_{H\in\mathcal{H}} {\frac{|V(H)|}{4}} = \frac{n}{4}.
\end{align*}
Suppose $\mathcal{H}' \neq \emptyset$. For every $x \in N(v)$ and $H \in \mathcal{H}$ such that $xy_{x,H} \in E(G)$ for some $y_{x,H} \in V(H)$, we say that \emph{$x$ is linked to $H$} and that \emph{$H$ is linked to $x$}. Since $G$ is connected, each member of $\mathcal{H}$ is linked to at least one member of $N(v)$. For each $x \in N(v)$, let $\mathcal{H}_{x} = \{H \in \mathcal{H} \colon H$ is linked to $ x\}$, $\mathcal{H}'_{x} = \{H \in \mathcal{H}' \colon H $ is linked to $x\}$ and $\mathcal{H}_{x}^* = \{H \in \mathcal{H}^* \colon H$ is linked to $x$ only$\}$. 

For each $H \in \mathcal{H}^*$, let $D_H$ be a smallest $P_3$-isolating set of $H$. For each $x \in N(v)$ and each $H \in \mathcal{H}_x'$, let $H_x' = H - y_{x,H}$, and let $D_{x,H}$ be a smallest $P_3$-isolating set of $H_x'$. By the induction hypothesis, $|D_{x,H}| \leq (|V(H)|-3)/4$.\medskip
\\
\textbf{Case 1:} $|\mathcal{H}'_x| \geq 2$ for some $x \in N(v)$. For each $H \in \mathcal{H}' \setminus \mathcal{H}'_x$, let $x_H \in N(v)$ such that $H$ is linked to $x_H$. Let $X = \{x_H \colon H \in \mathcal{H}' \setminus \mathcal{H}'_x\}$. Note that $x \notin X$. Let 
\[D = \{v, x\} \cup X \cup \left( \bigcup_{H \in \mathcal{H}_x'} D_{x,H} \right) \cup \left( \bigcup_{H \in \mathcal{H}' \setminus \mathcal{H}_x'} D_{x_H,H} \right) \cup \left( \bigcup_{H \in \mathcal{H}^*} D_{H} \right).\]
We have $V(G) = N[v] \cup \bigcup_{H \in \mathcal{H}} V(H)$, $y_{x,H} \in N[x]$ for each $H \in \mathcal{H}'_x$, and $y_{x_H,H} \in N[x_H]$ for each $H \in \mathcal{H}' \setminus \mathcal{H}'_x$, so $D$ is a $P_3$-isolating set of $G$. We have
\begin{align*} \iota(G, P_3) &\leq |D| = 2 + |X| + \sum_{H \in \mathcal{H}_x'} |D_{x,H}| + \sum_{H \in \mathcal{H}' \setminus \mathcal{H}_x'} |D_{x_H,H}| + \sum_{H \in \mathcal{H}^*} |D_{H}| \nonumber \\
&\leq \frac{8 + 4|X|}{4} + \sum_{H \in \mathcal{H}_x'} \frac{|V(H)| - 3}{4} + \sum_{H \in \mathcal{H}' \setminus \mathcal{H}_{x}'} \frac{|V(H)| - 3}{4} + \sum_{H \in \mathcal{H}^*} \frac{|V(H)|}{4} \nonumber \\ 
&= \frac{2+|X| - 3(|\mathcal{H}_{x}'| - 2) - 3(|\mathcal{H}'\setminus\mathcal{H}_{x}'| - |X|)}{4} + \sum_{H \in \mathcal{H}} \frac{|V(H)|}{4}.
\end{align*}
Thus, since $|\mathcal{H}_{x}'| \geq 2$, $|\mathcal{H}'\setminus\mathcal{H}_{x}'| \geq |X|$ and $2 + |X| = |\{v, x\} \cup X| \leq |N[v]|$, we have 
\[\iota(G,P_3) \leq \frac{2+|X|}{4} + \sum_{H\in\mathcal{H}} \frac{|V(H)|}{4} \leq \frac{|N[v]|}{4} + \sum_{H\in\mathcal{H}} \frac{|V(H)|}{4} = \frac{n}{4}.\]

\noindent\textbf{Case 2:} 
\begin{equation}|\mathcal{H}'_x| \leq 1 \mbox{ for each } x \in N(v). \label{beyondcase2} 
\end{equation}
For each $H \in \mathcal{H}'$, let $x_H \in N(v)$ such that $H$ is linked to $x_H$. Let $X = \{x_H \colon H \in \mathcal{H}'\}$. By (\ref{beyondcase2}), $|X| = |\mathcal{H}'|$. Let $W = N(v) \setminus X$ and
\[D = \{v\} \cup X \cup \left( \bigcup_{H \in \mathcal{H}'}D_{x_H,H} \right) \cup \left( \bigcup_{H \in \mathcal{H}^*} D_{H} \right).\] 
Then, $D$ is a $P_3$-isolating set of $G$.\medskip
\\
\textbf{Case 2.1:} $|W|\geq3$. Similarly to Case~1, we have
\begin{align*} \iota(G,P_3) &\leq |D| = 1 + |X| + \sum_{H \in \mathcal{H}'} |D_{x_H,H}| + \sum_{H \in \mathcal{H}^*} |D_{H}|  \\
&\leq \frac{|\{v\} \cup W| + 4|X|}{4} + \sum_{H \in \mathcal{H}'} \frac{|V(H)| - 3}{4} + \sum_{H \in \mathcal{H}^*} \frac{|V(H)|}{4} \\
&= \frac{|N[v]\setminus X| + |X|}{4} + \sum_{H \in \mathcal{H}} \frac{|V(H)|}{4}=\frac{n}{4}.
\end{align*} 

\noindent\textbf{Case 2.2:} $|W| \leq 2$.\medskip
\\
\textbf{Case 2.2.1:} Some member $H$ of $\mathcal{H}'$ is linked to $x_H$ only. Let $G^* = G - (\{x_H\} \cup V(H))$. Then, $G^*$ has a component $G_v^*$ such that $N[v] \setminus \{x_H\} \subseteq V(G_v^*)$, and any other component of $G^*$ is a member of $\mathcal{H}_{x_H}^*$. Let $D^*$ be a smallest $P_3$-isolating set of $G_v^*$. Then, $D^* \cup \{x_H\} \cup D_{x_H,H} \cup \bigcup_{I \in \mathcal{H}_{x_H}^*} D_I$ is a $P_3$-isolating set of $G$, so
\begin{align*} \iota(G, P_3) 
&\leq |D^*| + 1 + |D_{x_H,H}| + \sum_{I \in \mathcal{H}_{x_H}^*} |D_I| \\ 
&\leq \iota(G_v^*, P_3) + \frac{4}{4}+\frac{|V(H)|-3}{4} + \sum_{I \in \mathcal{H}_{x_H}^*}\frac{|V(I)|}{4} \\
&\leq \iota(G_v^*,P_3) + \frac{|\{x_H\}\cup V(H)|}{4} + \sum_{I \in \mathcal{H}_{x_H}^*} \frac{|V(I)|}{4}.
\end{align*}
This gives $\iota(G,P_3) \leq n/4$ if $\iota(G_v^*,P_3) \leq |V(G_v^*)|/4$. Suppose $\iota(G_v^*,P_3) > |V(G_v^*)|/4$. Then, $4 \iota(G_v^*,P_3) \geq |V(G_v^*)| + 1$, so $\iota(G_v^*,P_3) \geq (|V(G_v^*)| + 1)/4$. By the induction hypothesis, $\iota(G_v^*,P_3) = (|V(G_v^*)|+1)/4$. Let $G_v^- = G_v^* - v$. By the induction hypothesis, $\iota(G_v^-,P_3) \leq (|V(G_v^*)| - 3)/4$. Let $D^-$ be a smallest $P_3$-isolating set of $G_v^-$. Since $v \in N[x_H]$, $D^- \cup \{x_H\} \cup D_{x_H,H} \cup \bigcup_{I \in \mathcal{H}_{x_H}^*} D_I$ is a $P_3$-isolating set of $G$, so
\begin{align*} \iota(G, P_3) 
    &\leq |D^-| + 1 + |D_{x_H,H}| + \sum_{I \in \mathcal{H}_{x_H}^*} |D_I| \\ 
    &\leq \iota(G_v^-, P_3) + \frac{4}{4}+\frac{|V(H)|-3}{4} + \sum_{I \in \mathcal{H}_{x_H}^*}\frac{|V(I)|}{4} \\
    &\leq \frac{|V(G_v^*)|-3}{4} + \frac{|\{x_H\}\cup V(H)|}{4} + \sum_{I \in \mathcal{H}_{x_H}^*} \frac{|V(I)|}{4} < \frac{n}{4}.
\end{align*}
\noindent\textbf{Case 2.2.2:} For each $H \in \mathcal{H}'$, $H$ is linked to some $x_H' \in N(v) \setminus \{x_H\}$. Let $X' = \{x_H': H\in \mathcal{H}'\}$. By (\ref{beyondcase2}), for every $H, I \in \mathcal{H}'$ with $H \neq I$, we have $x_H' \neq x_I'$ and $x_H \neq x_I'$. Thus, $|X'| = |\mathcal{H}'|$ and $X' \subseteq W$. This yields $|\mathcal{H}'| \leq 2$ as $|W| \leq 2$.

Suppose $|\mathcal{H}'| = 2$. Then, $|W| = 2$ and $X' = W$. Let $H_1$ and $H_2$ be the two members of $\mathcal{H}'$. We have $N(v) = \{x_{H_1},x'_{H_1},x_{H_2},x'_{H_2}\}$. Thus, $\Delta(G) = 4$. Therefore, by Lemma~\ref{LemIsol1}, Lemma~\ref{LemIsol2} and the induction hypothesis,
\begin{align*}
    \iota(G,P_3) &\leq 1 + \iota(H_1,P_3) + \iota(H_2,P_3) + \sum_{H\in\mathcal{H}^*}\iota(H,P_3) \\
    &\leq \frac{|N[v]|-1}{4} + \frac{|V(H_1)| + 1}{4} + \frac{|V(H_2)| + 1}{4} + \sum_{H\in\mathcal{H}^*}\frac{|V(H)|}{4} = \frac{n+1}{4}.
\end{align*}

Now suppose $|\mathcal{H}'| = 1$. Let $H_1$ be the member of $\mathcal{H}'$. Since $3 \leq \Delta(G) = d(v) = |X| + |W| = 1 + |W| \leq 3$, $\Delta(G) = 3$. Therefore, by Lemma~\ref{LemIsol1}, Lemma~\ref{LemIsol2} and the induction hypothesis,
\begin{align*}\pushQED{\qed}
    \iota(G,P_3) &\leq 1 + \iota(H_1,P_3) + \sum_{H\in\mathcal{H}^*}\iota(H,P_3) \\
    &\leq \frac{|N[v]|}{4} + \frac{|V(H_1)| + 1}{4} + \sum_{H\in\mathcal{H}^*}\frac{|V(H)|}{4} = \frac{n+1}{4}. \qedhere
    \popQED
\end{align*}
\\

\noindent \textbf{Acknowledgements.} The authors are grateful to the anonymous referee for checking the paper carefully and providing constructive remarks. The research work disclosed in this publication is partially funded by the Tertiary Education Scholarships Scheme (Malta).

\end{document}